\documentclass[12pt,a4paper]{article}
\usepackage{amsfonts, latexsym, amsthm, color}

\def\seq#1#2#3{#1_{#2},\,\ldots,#1_{#3}}

\def\VV{{\underline{V}}}

\def\vv{{\underline{v}}}
\def\tt{{\underline{t}}}

\def\mm{\underline{m}}

\def\kk{\underline{k}}
\def\1{\underline{1}}

\def\P{\Bbb P}

\def\Z{\Bbb Z}
\def\Q{\Bbb Q}
\def\C{\Bbb C}

\def\OO{{\cal O}}

\newtheorem{theorem}{Theorem}

\newtheorem{proposition}{Proposition}

\newenvironment{definition}
{\smallskip\noindent{\bf Definition\/}:}{\smallskip\par}

\newenvironment{examples}
{\smallskip\noindent{\bf Examples\/}.}{\smallskip\par}
\newenvironment{remark}
{\smallskip\noindent{\bf Remark\/}.}{\smallskip\par}
\newenvironment{remarks}
{\smallskip\noindent{\bf Remarks\/}.}{\smallskip\par}

\title{Poincar\'e series of collections  of plane valuations
\footnote{
Partially supported by the grant MTM2007-64704. Third author is
also partially supported by the grants RFBR-007-00593 and NSh-709.2008.1.
}
}
\author{
A.~Campillo
\and F.~Delgado
\and S.M.~Gusein-Zade
\and F. Hernando
\thanks{Addresses:
A. Campillo and F. Delgado:
University of Valladolid,
Dept. of Algebra, Geometry and Topology, 47011 Valladolid, Spain. E-mail:
campillo\symbol{'100}agt.uva.es, fdelgado\symbol{'100}agt.uva.es;
S.M.Gusein-Zade: Moscow State University, Faculty of
Mathematics and Mechanics, Moscow, GSP-1, 119991, Russia. E-mail:
sabir\symbol{'100}mccme.ru;
F. Hernando:
Dept. of Mathematics. University College Cork. Ireland.
E-mail: f.hernando\symbol{'100}ucc.ie
}}

\date{}
\begin{document}
\sloppy
\def\eps{\varepsilon}

\maketitle
\begin{abstract}
In earlier papers there were given formulae for the Poincar\'e
series of multi-index
filtrations on the ring $\OO_{\C^2,0}$ of germs of functions of
two variables defined by collections of valuations corresponding
to (reducible) plane curve singularities and by collections of
divisorial
ones. It was shown that the Poincar\'e series of a collection of
divisorial valuations determines the topology of the collection
of divisors.
Here we give a formula for
the Poincar\'e series of a general collection of valuations on
the ring $\OO_{\C^2,0}$ centred at the origin and prove
a generalization of the statement that the Poincar\'e series
determines the topology of the collection.
\end{abstract}

In \cite{duke}, \cite{divisorial}, ... there were considered
multi-index filtrations on the ring $\OO_{\C^2,0}$ of germs of
functions of two variables defined by collections of valuations
corresponding to (reducible) plane curve singularities and by
collections of divisorial ones. One gave formulae for the Poincar\'e
series of such filtrations in terms of an embedded resolution of the
curve singularity or of the collections of divisors respectively.
These formulae give Poincar\'e series as rational functions equal to
products/ratios of cyclotomic polynomials. In particular, it was
shown that the Poincar\'e series of the collection of valuations
corresponding to a curve coincides with the Alexander polynomial of
the corresponding algebraic link. This implies that the Poincar\'e
series of such a collection determines the topology of the curve
(\cite{yamamoto}). An analogue of this statement for divisorial
valuations was proved in \cite{FAOM}. Here we give a formula for the
Poincar\'e series of a general collection of valuations on the ring
$\OO_{\C^2,0}$ centred at the origin. We also prove a generalization
of the statement that the Poincar\'e series determines the topology
of the collection of divisors.

\section{Poincar\'e series of several valuations}
Let ${\cal O}_{X,0}$ be the ring of germs of functions on a
germ $(X,0)$ of a complex analytic variety and let $ G$ be
a ordered abelian group, $ G^+:=\{ a\in G:
 a\ge 0\}$. A valuation on the ring ${\cal O}_{X,0}$ with
values in the group $ G$ is a map $v: {\cal O}_{X,0}\to
 G^+\cup\{+\infty\}$ such that:
\begin{enumerate}
\item $v(g_1\cdot g_2)=v(g_1)+v(g_2)$;
\item $v(g_1+g_2) \ge \min\{v(g_1),v(g_2)\}$;
\item $v(\lambda)=0$ \ for \ $\lambda\in\C^*=\C\setminus\{0\}$.
\end{enumerate}
The semigroup $S=\mbox{Im\,}v\setminus\{+\infty\}$ (the
semigroup of values of the valuation $v$) is well ordered (i.e.
each subset of $S$ has the minimal element) and moreover each
element $a\in S$ has a finite number of representations as the
sum $a_1+a_2$ of two elements of $S$ (\cite{zariski}, see also
\cite{CG}).

Let $v_i$, $i=1, \ldots, s$, be valuations on the ring ${\cal
O}_{X,0}$ with values in ordered groups $ G_i$ and the
semigroups of values $S_i$. The direct product
$S=S_1\times\ldots\times S_s$ is a partially ordered semigroup:
$\vv'=(v_1', \ldots,v_s')\ge \vv''=(v_1'', \ldots,v_s'')$
iff $v_i' \ge v_i''$ for all $i=1, \ldots, s$ ($v_i',
v_i''\in S_i$). For a germ $g\in {\cal O}_{X,0}$, let
$\vv(g):=(v_1(g), \ldots, v_s(g))$.

\begin{definition}
The ring $\Z[[S]]$ of power  series on the semigroup
$S=S_1\times\ldots\times S_s$  is the set of formal expressions
of the form $\sum\limits_{\vv\in S}a_{\vv}\, \tt^{\,\vv}$
($\vv:=(v_1, \ldots, v_s) \in S$, $\tt^{\,\vv}
:= t_1^{\,v_1}\cdot\ldots\cdot
t_s^{\,v_s}$) with the usual ring operations.
\end{definition}

\begin{remarks}
{\bf 1.} The fact that $\Z[[S]]$ is a ring (i.e. that the
multiplication is defined) follows from the described
properties of the semigroup of values of a valuation.

\noindent {\bf 2.} If $ G_i=\Z$ for $i=1, \ldots, s$, the
ring $\Z[[S]]$ is contained  in the ring $\Z[[t_1, \ldots,
t_s]]$ of power series in the variables $t_1, \ldots, t_s$ with
integer coefficients and thus an element of the ring $\Z[[S]]$
is a usual power series in the variables $t_1, \ldots, t_s$.
\end{remarks}

The collection $\{v_i\}$ of valuations defines a multi-index
filtration on the ring ${\cal O}_{X,0}$ by the ideals
$$
J(\vv)= \{g\in {\cal O}_{X,0}: \vv(g)\ge\vv\}
$$
(indexed by the elements $\vv$ of the group
$ G= G_1\times\ldots\times G_s$). For $I\subset
I_0=\{1, \ldots, s\}$, $\vv=(v_1, \ldots, v_s)\in G$, let
$$
J^{+I}(\vv) := \{g\in J(\vv): v_i(g)>v_i \mbox{\  for \ }i\in
I\}\,,
$$
$ J^{+}(\vv) := J^{+I_0}(\vv)$.

\begin{definition}
The Poincar\'e series $P_{\{v_i\}}(t_1, \ldots, t_s)$ of the
collection of valuations $\{v_i\}$ is the element of the ring
$\Z[[S]]$ defined by
$$
P_{\{v_i\}}(\tt)= \sum\limits_{\vv\in S}\left(\sum\limits_{I\subset
I_0}(-1)^{\#I}\dim\left(J^{+I}(\vv)/J^{+}(\vv)\right)\right)\cdot
\tt^{\,\vv}\,.
$$
\end{definition}

One can see that, for collections of integer valued valuations,
this definition coincides with that used e.g. in \cite{duke}. One
can easily extract a proof of this from the proof of Theorem~3 in
\cite{duke}.

\begin{remark}
This notion is defined if all the factor spaces
$J^{+I}(\vv)/J^{+}(\vv)$ have finite dimensions.
This takes place if each valuation $v_i$, $i=1,\ldots, s$, is
centred at the origin, i.e. $\{g\in {\cal O}_{X,0} :
v_i(g)>0\}$ coincides with the maximal ideal ${\frak m}$ of the
ring ${\cal O}_{X,0}$.
\end{remark}

\begin{definition}
A valuation $v$ with values in a group $ G$ is finitely
determined if, for each $v_0\in G$, the condition
$v(g)=v_0$ is a constructible condition on a jet of the germ
$g$ of a certain (finite) order (see \cite{Int}).
\end{definition}

\begin{examples}
{\bf 1.} For an irreducible plane curve germ $C=\{f=0\}\subset
(\C^2,0)$, let $\varphi: (\C,0) \to (\C^2,0)$ be a
parametrization (uniformization) of the curve $C$. For a germ
$g\in {\cal O}_{X,0}$, let $v_C(g)$ be the order of zero of the
function $g\circ\varphi(\tau)$ at the origin (if
$g\circ\varphi \equiv 0$, $v_C(g):=+\infty$). The map $v_C: {\cal
O}_{X,0}\to\Z\cup\{+\infty\}$ is a finitely determined rank $1$
valuation on the ring ${\cal O}_{X,0}$. (This valuation has a
non-trivial kernel, i.e. $v_C(g)=+\infty$ for some $g\ne 0$.)

\noindent {\bf 2.} For the same curve $C=\{f=0\}$ and for $g\in
{\cal O}_{X,0}$, $g\ne 0$, let $g=f^{k_C(g)} g'$, where $g'$
is not divisible by $f$. The map $w_C: {\cal
O}_{X,0}\to\Z^2\cup\{+\infty\}$, $w_C(g):=(k_C(g), v_C(g'))$, is
a rank $2$ valuation on the ring ${\cal O}_{X,0}$. (The group
$\Z^2$ is ordered lexicographically, i.e. $(k',v')>(k'',v'')$
iff $k' >k''$ or $k' =k''$ and $v'>v''$.) The valuation $w_C$
is not finitely determined: for $k>0$, the condition $k_C(g)=k$
is not determined by a (finite) jet of the germ $g$. (The
kernel of the valuation $w_C$ is trivial.)
\end{examples}

The notion of integration with respect to the Euler
characteristic $\chi$ over the projectivization $\P{\cal
O}_{X,0}$ of the ring ${\cal O}_{X,0}$ was defined in \cite{Int}.
The argument from \cite{Int} give the following statement.

\begin{proposition}
If all the valuations $v_i$, $i=1, \ldots, s$, are finitely
determined, one has
\begin{equation}\label{int}
P_{\{v_i\}}(\tt)=\int_{\P{\cal O}_{X,0}}\tt^{\,\vv(g)}d\chi \,.
\end{equation}
\end{proposition}

\begin{remark}
One can formulate the notion of the Poincar\'e series of
a collection of valuations in terms of the extended semigroup of
the collection in spirit of \cite{duke}, \cite{divisorial}.
\end{remark}

\section{Valuations on ${\cal O}_{\C^2,0}$}

Let $\pi: (X,D)\to(\C^2,0)$ be
a modification of the plane
by a (finite) sequence of blowing-ups. The exceptional divisor
$D$ is the union of irreducible components $E_{\sigma}$
($\sigma\in\Gamma$); each of them is isomorphic to the complex
projective line $\C\P^1$. The dual graph of the modification
$\pi$ consists of vertices corresponding to the irreducible
components $E_{\sigma}$, i.e. to elements of $\Gamma$; two
vertices are connected
by an edge iff the corresponding components intersect. To
$\sigma\in\Gamma$, i.e. to an irreducible component
$E_{\sigma}$ of the exceptional divisor, there corresponds a
natural valuation $v_{\sigma}$: the {\em divisorial} one. For a
germ
$g\in {\cal O}_{X,0}$ the value $v_{\sigma}(g)$ is defined as a
multiple of (say, $c_{\sigma}$ times) the order $w_{\sigma}(g)$
of zero of the function $g\circ\pi$ along the component
$E_{\sigma}$ (i.e. of the coefficient at $[E_{\sigma}]$ in the
zero divisor of the function $g\circ\pi$). It is convenient to
choose the coefficient $c_{\sigma}$ in such a way that $\min\limits_
{g\in{\frak m}} v_{\sigma}(g)=1$, i.e.
$c_{\sigma}=1/w_{\sigma}(g)$ for a generic function $g\in{\frak
m}$.

It is known (\cite{zariski}, see also \cite{DGN}) that all
valuations on the ring
${\cal O}_{X,0}$ centred at the origin (i.e. such that $v(g)>0$ for
$g\in{\frak m}$) are given by the following list.

{\bf I.} All rank $1$ valuations centred at the origin
correspond to some sequences (finite or infinite) of blowing-ups
such that each
next blowing-up is made at a point of the divisor born on the
previous step. To get a one-to-one correspondence, one should
exclude sequences of blowing-ups made at each step
after a certain one at the intersection point of a fixed
divisor $E_{\sigma_0}$ with the last one. (In such case
the correspondence described below leads to the divisorial
valuation $v_{\sigma_0}$ associated to the divisor
$E_{\sigma_0}$.) If the sequence is finite, the corresponding
valuation is the divisorial one associated to the last divisor.
If the sequence is infinite, the value $v(g)$ of the
corresponding valuation
is defined as the limit $\lim\limits_{i\to\infty}
v_{\sigma_i}(g)$ where
$v_{\sigma_i}$ is the divisorial valuation
associated to the divisor born on $i$-th step. Depending on the
sequence of blowing-ups, one can distinguish the following types of valuations.

{\bf I.1.} The blowing-ups are made at the intersection points
of the strict transforms of a fixed irreducible curve
$C=\{f=0\}\subset (\C^2, 0)$ with the exceptional divisor. The
corresponding valuation
(a
{\em curve valuation of rank 1}) is equivalent to the valuation
defined by the order of a function $g$ on the curve $C$ in an
uniformization parameter (Example 1 above). The corresponding
dual graph of the modification has a growing infinite tail:
Fig.1. This valuation has a non-trivial kernel.

\begin{figure}[h]
$$
\unitlength=1.00mm
\begin{picture}(80.00,20.00)(-10,13)
\thicklines \put(-5,30){\line(1,0){41}}
\put(44,30){\line(1,0){12}} \put(38,30){\circle*{0.5}}
\put(40,30){\circle*{0.5}} \put(42,30){\circle*{0.5}}
\put(30,10){\line(0,1){20}} 
\put(10,15){\line(0,1){15}}
\thinlines
\put(20,30){\circle*{1}} \put(30,30){\circle*{1}}
\put(50,30){\circle*{1}}
\put(58,30){\circle*{0.5}}
\put(60,30){\circle*{0.5}}
\put(62,30){\circle*{0.5}}
\put(30,20){\circle*{1}}
\put(10,30){\circle*{1}}
\put(30,10){\circle*{1}}
\put(-5,30){\circle*{1}}
\put(0,30){\circle*{1}} \put(5,30){\circle*{1}}
\put(15,30){\circle*{1}} \put(25,30){\circle*{1}}
\put(35,30){\circle*{1}} \put(45,30){\circle*{1}}
\put(55,30){\circle*{1}} \put(10,25){\circle*{1}}
\put(10,20){\circle*{1}} \put(10,15){\circle*{1}}
\put(30,25){\circle*{1}} \put(30,15){\circle*{1}}
\put(35,30){\circle*{1}}
\put(-9,29){{\bf 1}}
\end{picture}
$$
\caption{The dual graph of a valuation of type I.1.}
\label{fig1}
\end{figure}
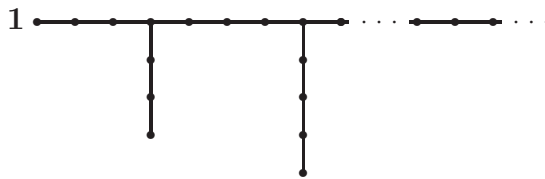

{\bf I.2.} The sequence of blowing-ups produces an infinite tail
like on Fig.1, but there exists no curve corresponding to this
sequence (one can say that there is a formal curve defined by a
formal power series), one gets a discrete valuation with the
trivial kernel: a {\em formal curve valuation}.

{\bf I.3.} The sequence of blowing-ups increases the number of
rupture points producing, as the limit, the graph shown on Fig.2.
In this case the group of values is contained in the ring $\Q$
of rational numbers and is not finitely generated. For any
$g\in \OO_{\C^2,0}$, one has $v(g)=v_{\sigma_i}(g)$ for $i$
sufficiently
large. Valuations of this sort we call {\em infinite valuations}.

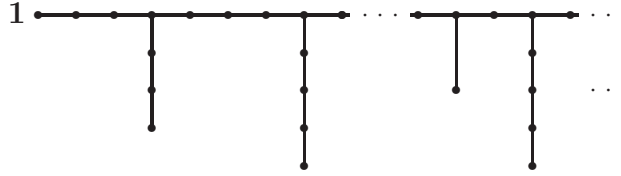
\begin{figure}[h]
$$
\unitlength=1.00mm
\begin{picture}(80.00,20.00)(-10,13)
\thicklines \put(-5,30){\line(1,0){41}}
\put(44,30){\line(1,0){22}} \put(38,30){\circle*{0.5}}
\put(40,30){\circle*{0.5}} \put(42,30){\circle*{0.5}}
\put(30,10){\line(0,1){20}} \put(50,20){\line(0,1){10}}
\put(60,10){\line(0,1){20}} \put(10,15){\line(0,1){15}}
\thinlines
\put(20,30){\circle*{1}} \put(30,30){\circle*{1}}
\put(50,30){\circle*{1}} \put(60,30){\circle*{1}}
\put(65,30){\circle*{1}}
\put(68,30){\circle*{0.5}}
\put(70,30){\circle*{0.5}}
\put(72,30){\circle*{0.5}}
\put(68,20){\circle*{0.5}}
\put(70,20){\circle*{0.5}}
\put(72,20){\circle*{0.5}}
\put(30,20){\circle*{1}} \put(60,25){\circle*{1}}
\put(60,20){\circle*{1}} \put(60,15){\circle*{1}}
\put(10,30){\circle*{1}} \put(30,10){\circle*{1}}
\put(50,20){\circle*{1}} \put(60,10){\circle*{1}}
\put(-5,30){\circle*{1}} \put(0,30){\circle*{1}}
\put(5,30){\circle*{1}} \put(15,30){\circle*{1}}
\put(25,30){\circle*{1}} \put(35,30){\circle*{1}}
\put(45,30){\circle*{1}} \put(55,30){\circle*{1}}
\put(10,25){\circle*{1}} \put(10,20){\circle*{1}}
\put(10,15){\circle*{1}} \put(30,25){\circle*{1}}
\put(30,15){\circle*{1}} \put(35,30){\circle*{1}} \put(-9,29){{\bf
1}}
\end{picture}
$$
\caption{The dual graph of a valuation of type I.3.} \label{fig2}
\end{figure}

{\bf I.4.} In the cases I.2 and I.3 blowing-ups were made either
permanently after a certain step or from time to time at a smooth
point of the last divisor. If this is not the case (i.e. if all
the blowing-ups starting from a  moment are produced at an
intersection point of the last divisor with a previous one but
not with one and the same) one gets a valuation with a
finitely generated group of values, which is not contained in
$\Q$ (an {\em irrational valuation}).

{\bf II.} Besides that there are two kinds of rank 2 valuations
on the ring $\OO_{\C^2,0}$.

{\bf II.1.} For an irreducible curve germ
$C=\{f=0\}\subset (\C^2,0)$ one has the valuation
described in Example 2: {\em curve valuations of rank 2}.
Valuations
of this type are the only ones on the list which are not finitely
determined.

{\bf II.2.} For a divisor $E_{\sigma}$ of a modification
$\pi: (X,D)\to (\C^2,0)$ and for a point $P$ on it, let
$E_{\sigma}$ be locally given by an equation $u=0$. Then we can
define a rank $2$ valuation $v = v(g)$ with values in the group
$\Z^2$ ordered lexicographically applying the construction
of II.1 to the germ of function $g\circ \pi$ at the point $P$ of
the curve $u=0$. This means that one writes $g\circ \pi$ as
$u^{v_{\sigma}(g)} g'$ where $g'\notin \langle u\rangle$
($v_{\sigma}(g)$ is just the value of $g$ with respect to the
divisorial valuation $v_{\sigma}$)  and puts
$v(g) := (v_{\sigma}(g), w_{u}(g'))$, where $w_u(g')$ is the
order of zero of the function $g'$ on the smooth curve $u=0$ at
the point $P$. This valuation is called an {\em exceptional curve
valuation}.

If one excludes valuations of type I.1. (i.e. those which have
a non-trivial kernel), one can say that all valuations on the
ring $\OO_{\C^2,0}$ are in one to one correspondence with all
(finite and infinite) sequences of blowing-ups of the described
type. In this case sequences of blowing-ups excluded in I
correspond to valuations of type II.2. To a  sequence of the
type described in I.1 one associates the valuation II.1 (of
rank 2) corresponding to the curve germ.

\begin{remarks}
{\bf 1.} A valuation of type I.1 is essentially the second
component of the corresponding valuation of type II.1. If one
considers the ring $\C[[x,y]]$ of formal power series in two
variables instead of $\OO_{\C^2,0}$, valuations of type I.1 and
I.2 constitute one and the same type. The results of the paper
are valid in this setting as well.

\noindent{\bf 2.} Except valuations centred at the origin, there
are so called $f$-adic valuations $k_C$ corresponding to
irreducible plane curve singularities $C=\{f=0\}\subset
(\C^2,0)$. For a germ $g\in \OO_{\C^2,0}\setminus\{0\}$,
$k_C(g)$
is defined by the relation $g= f^{k_{C}(g)} g'$ where $g'\notin
\langle f\rangle$. This is just the first component of the rank
2 curve valuation corresponding to $C$. For an appropriate
definition, the Poincar\'e series of a collection of valuations
containing $f$-adic ones is a reduction of the Poincar\'e series
of the collection obtained by substituting $f$-adic ones by the
corresponding plane curve valuations of rank 2.
\end{remarks}

For further discussions it is convenient to use a notion of a
resolution of  a (finite) collection of valuations. For a divisorial
valuation this is a modification $\pi:(X,D)\to(\C^2,0)$ of the plane
by a finite sequence of blowing-ups over the origin which contains
the corresponding divisor (this means that if is a modification of
the minimal resolution by a finite sequence of blowing-ups). For a
valuation of one of the types I, II.1, let $\pi_i: (X_i,D_i)\to
(\C^2,0)$ be the modification obtained at the $i$-th step of the
corresponding sequence of blowing-ups. Let $(X,D)$ be the projective
limit (in the category of analytical spaces) of the sequence
$\{(X_i,D_i)\}$ and let $\pi: (X,D)\to (\C^2,0)$ be the
corresponding map. One can see that $X$ is a smooth complex surface
and $D$ is the union of infinitely many projective lines on it. The
map $\pi: (X,D)\to (\C^2,0)$ will be called the minimal resolution
of the valuation. This is not a resolution since since, in
particular, the map $\pi$ is not proper. For a valuation of type
II.2, the minimal resolution will be defined as the minimal
resolution of the corresponding divisorial valuation (the first
component of the considered one) followed by an additional
blowing-up at the corresponding point. Thus in the minimal
resolution the point under consideration is an intersection point of
components of the exceptional divisor.

The minimal resolution of a finite
collection of valuations is the projective limit of the
corresponding (multi-index) system of modifications. It is simply
the fibre product of the minimal resolutions of all valuations.

Finally a resolution of a finite collection of valuations is a
modification of the space $(X,D)$ of the minimal resolution by a
finite number of blowing-ups (at points of $D$).

\section{Poincar\'e series of a collection of plane valuations}

Let $v_i$, $i=1,\ldots, s$, be a set of valuations on the ring
$\OO_{\C^2,0}$ centred at the origin. Let $\pi:(X,D)\to(\C^2,0)$
be a resolution of the
set of valuations $\{v_i\}$. For a component $E_\sigma$ of the
exceptional divisor $D$, $\sigma\in \Gamma$, let $L_\sigma$ be a
smooth germ of curve
transversal to $D$ at a smooth point of it, let the curve germ
$\ell_\sigma = \pi(L_\sigma)\subset (\C^2,0)$ be given by an
equation $\{g_\sigma=0\}$, and let
$m^{\sigma}_{i} =v_i(g_\sigma)$,
$\mm^{\sigma} := (\seq{m^{\sigma}}1s)\in  G$.

Without loss of generality one can suppose that $\seq{v}1r$ are
valuations of type II.1 (or rank 2) and the others are  not. For
a component $E_\sigma$ of the exceptional divisor $D$, let
$\stackrel{\circ}{E_\sigma}$ be its ``smooth part", i.e.
$E_\sigma$ itself minus the intersection points with all other
components of $D$. For $i=1,\ldots,r$, let $C_i=\{f_i=0\}$ be the
corresponding curve, $f_i\in \OO_{\C^2,0}$.

\begin{theorem}
The Poincar\'e series of the set of valuations $\{v_i\}$ is
given by the equation
\begin{equation}\label{eq2}
P_{\{v_i\}}(\tt) = \prod_{\sigma\in \Gamma}
(1-\tt^{\,\mm^\sigma})^{-\chi(\stackrel{\circ}{E_\sigma})}
\times
\prod_{i=1}^r (1- t_i^{(1,0)}\prod\limits_{j\neq
i}t_j^{v_j(f_i)})^{-1}
\; .
\end{equation}
\end{theorem}

\begin{proof}
One cannot use the equation~(\ref{int}) directly since valuations
of type II.1 are not finitely determined. However the condition
$v_i(g)=(0,v)$ is defined by a jet of the germ $g$ of a finite
order. Therefore let  us first compute the part $P^{o}(\tt)$ of
the Poincar\'e series $P_{\{v_i\}}(\tt)$ which consists of
monomials which are not divisible by any of $t_i^{(1,0)}$,
$i=1,\ldots,r$. One can see that
$$
P^o(\tt) = \int\limits_{\P(\OO_{\C^2,0}\setminus
\bigcup\limits_{i=1}^r \langle f_i\rangle)} \tt^{\,\vv} d \chi
$$

Computation of $P^o(\tt)$ essentially repeats the
arguments from, e.g., \cite{Int} or
\cite{divisorial}. To compute the series $P^o(\tt)$ up to
terms of any
fixed degree $\VV\in  G$ one can make finitely many additional
blowing-ups at intersection points of the components of the
exceptional divisor $D$ so that, for a germ $g\in
\OO_{\C^2,0}\setminus \bigcup\limits_{i=1}^r \langle f_i\rangle$ with
$\vv(g)\le \VV$, all the intersection points of the
strict transform of the curve $\{g=0\}$ with the exceptional
divisor $D$ belong to its smooth part
$\stackrel{\circ}{D}=\bigcup\limits_{\sigma}\stackrel{\circ}{E_\sigma}$.
Let $\OO^{\VV}_{\C^2,0} = \{g\notin \bigcup\limits_{i=1}^r
\langle f_i\rangle :  \vv(g)\le \VV\}$, and let
$\OO^{D}_{\C^2,0}$ be the set of germs
$g\in \OO_{\C^2,0}\setminus \{0\}$ such that  the intersection
points of the strict transform of the curve $\{g=0\}$ with the
exceptional divisor $D$ belong to
$\stackrel{\circ}{D}$. The integral
$$
\int_{\P\OO^{D}_{\C^2,0}} \tt^{\,\vv} d \chi
$$
over the projectivization $\P\OO^{D}_{\C^2,0}$ contains all terms of
the series $P^o(\tt)$ up to degree $\VV$. There is a map from
$\P\OO^{D}_{\C^2,0}$ to the space of effective divisors on
$\stackrel{\circ}{D}$: to a function germ $g\in \OO^D_{\C^2,0}$ one
associates the intersection of the strict transform of the curve
$\{g=0\}$ with the exceptional divisor $D$. Proposition 2 from
\cite{Int} implies that the preimage of a point with respect to this
map is a complex affine space and thus has the Euler characteristic
equal to 1. The Fubini formula implies that the integral
$\int_{\P\OO^{D}_{\C^2,0}} \tt^{\,\vv} d \chi$ is equal to the
integral with respect to the Euler characteristic of the monomial
$\tt^{\,\vv}$ over the space of effective divisors on
$\stackrel{\circ}{D}$, where $\vv$ is the additive function on the
space of effective divisors (with values in $ G$) equal to
$\mm^\sigma$ for a point  from the component
$\stackrel{\circ}{E_\sigma}$ of $\stackrel{\circ}{D}$.

The space of effective divisors on $\stackrel{\circ}{D}$ is the
direct product of the spaces of effective divisors on the
components $\stackrel{\circ}{E_\sigma}$, $\sigma\in \Gamma$. Each
of the latter ones is the disjoint union of the symmetric powers
$S^k\stackrel{\circ}{E_\sigma}$ of the component
$\stackrel{\circ}{E_\sigma}$. Therefore
$$
\int_{\P\OO^{D}_{\C^2,0}} \tt^{\,\vv} d \chi = \prod_{\sigma\in
\Gamma}
\left( \sum_{k=0}^{\infty} \chi(S^k \stackrel{\circ}{E_\sigma})
\cdot \tt^{\,k \mm^{\sigma}} \right) \; .
$$
Using the equation
$$
\sum_{k=0}^\infty  \chi(S^k X)\cdot t^k = (1-t)^{\chi(X)}
$$
one gets
\begin{equation}\label{eq3}
\int_{\P\OO^{D}_{\C^2,0}} \tt^{\,\vv} \; d \chi = \prod_{\sigma}
(1-\tt^{\,\mm^\sigma})^{-\chi(\stackrel{\circ}{E_\sigma})}\; .
\end{equation}
The right hand side of equation~\ref{eq3} do not contain
components of the exceptional divisor born under additional
blowing-ups since for each of them
$\chi(\stackrel{\circ}{E_\sigma}) =0$. Therefore
$$
P^o (\tt) = \prod_{\sigma\in \Gamma}
(1-\tt^{\,\mm^\sigma})^{-\chi(\stackrel{\circ}{E_\sigma})}\;
$$
(for a resolution $\pi: (X,D)\to (\C^2,0)$ of the set of
valuations $\{v_i\}$).

Now, for $\kk = (\seq k1r)\in \Z^{r}_{\ge 0}$, let
$P^{\,\kk}(\tt)$ be the sum of terms of the  Poincar\'e series
$P_{\{v_i\}}(\tt)$ divisible by
$\prod\limits_{i=1}^r t_i^{(k_i,0)}$ but not divisible by a
monomial of this sort of higher degree. The set of functions $g$
with $\vv(g)$ divisible by
$\prod\limits_{i=1}^r t_i^{(k_i,0)}$ (but not by a monomial of
higher degree) is just
$\prod\limits_{i=1}^r f_i^{k_i} \cdot \left(
\OO_{\C^2,0}\setminus
\bigcup\limits_{i=1}^r \langle f_i\rangle\right)$. One has
$v_j(\prod\limits_{i=1}^r f_i^{k_i}) = \sum\limits_{i=1}^r k_i
v_j(f_i)$ and $v_i(f_i) = (1,0)$ for $i=1,\ldots, r$.

From this it follows that
$$
P^{\,\kk}(\tt) =
P^o (\tt)
\prod_{i=1}^r t_i^{(k_i,0)} \prod_{j\neq i} t_j^{\,k_i v_j(f_i)}
$$
and therefore
$$
P_{\{v_i\}}(\tt) = P^o (\tt) \sum_{\kk\in \Z^r_{\ge 0}}
\prod_{i=1}^r t_i^{(k_i,0)} \prod_{j\neq i} t_j^{\,k_i v_j(f_i)}
=
P^o(\tt) \cdot  \prod_{i=1}^r (1- t_i^{(1,0)}\prod\limits_{j\neq
i}t_j^{\,v_j(f_i)})^{-1}  \; .
$$
\end{proof}

\section{Poincar\'e series determines dual graphs}

To a resolution $\pi: (X,D)\to (\C^2,0)$ of a collection of
valuations one associates the dual graph (generally speaking
infinite). It consists of vertices corresponding to the irreducible
components of the exceptional divisor $D$; two vertices are
connected by an edge if the corresponding components intersect. The
set of vertices of the dual graph inherits a partial order defined
by approximation of the modification by sequences of blowing-ups: a
component $E_{\sigma'}$ is ``greater" than another component
$E_\sigma$ ($\sigma'>\sigma$) if the exceptional divisor of the
minimal modification which contains $E_{\sigma'}$ also contains
$E_\sigma$. Two resolutions are called combinatorially equivalent if
their dual graphs together with the partial orders are isomorphic.

For collections of valuations of types I.1 and I.2 the
Poincar\'e series
coincides with the Alexander polynomial (in several variables)
of the corresponding
algebraic link (obtained by cutting the non-convergent series
for a valuation of type I.2) (see \cite{duke}). It is known that
the
Alexander polynomial of an algebraic link determines the topology
of the curve singularity (\cite{yamamoto}, see another proof in
\cite{FAOM}). Therefore the Poincar\'e
series of a collection of such valuations determines its minimal
resolution up to combinatorial equivalence.

In \cite{FAOM}, it was shown that the Poincar\'e series of a
collection
of divisorial valuations ``determines the topology of the set of
divisors" in the sense that it determines the dual graph of the
minimal resolution up to combinatorial equivalence. Moreover,
it was shown that the Poincar\'e series of a collection of
valuations which includes both divisorial ones and those of types
I.1 (or I.2) does not determine, in general, the dual graph of
the minimal resolution.

\begin{theorem}
Suppose that a collection $\{v_i\}$
of valuations does not contain
valuations if types I.1 and I.2. Then the Poincar\'e series
$P_{\{v_i\}}(\tt)$ of the collection determines types of the
valuations, the
dual graph of the minimal resolution up to combinatorial equivalence and divisors or
sequences of divisors corresponding to the valuations.
\end{theorem}

In what follows, for short, we shall simply say that the
Poincar\'e series determines the dual graph of the minimal
resolution.

\begin{proof}
One property of collections of valuations of the type
under consideration used in the proof is the fact that the
projection formula holds for them. This means that, for a
subset $I\subset I_0=\{1, \ldots, s\}$, the Poincar\'e series
of the set of valuations $\{v_i\}_{i\in I}$ is obtained from
the Poincar\'e series of the whole set $\{v_i\}$ by omitting
the variables $t_i$ with $i\notin I$ (in other words by
substituting $t_i$ with $i\notin I$ by $1$).

\begin{remark}
Omitting a valuation of type I.3 (an infinite valuation), as a
factor one formally gets an infinite product of the form
$$
(1-\tt^{\,\mm^1})(1-\tt^{\,\mm^1})^{-1}(1-\tt^{\,\mm^2})
(1-\tt^{\,\mm^2})^{-1}\ldots
$$
This product should be canceled.
\end{remark}

The projection formula follows directly from the equation
(\ref{eq2}) (Theorem 1).

\begin{remark}
For collections of valuations which do
not contain valuations of type II.1 (curve valuations of rank 2)
this can be also deduced from equation (\ref{int}). The fact that
the projection formula holds for valuations of type II.1 as well
can mean that this valuations can be considered finitely
determined in some weak sense.
\end{remark}

The projection formula implies that, from the Poincar\'e series
 of a collection of valuations, one can restore, in particular,
the Poincar\'e series of each individual valuation from the
collection.

For each individual valuation one can define its type
from the Poincar\'e series. From equation~(\ref{eq2}) it follows
that the set of exponents in the Poincar\'e series with
non-vanishing coefficients generates an abelian group of rank
equal to the rank of the valuation. The Poincar\'e series of a
valuation of type II.1 (a curve valuation of rank 2) has
infinitely many non-vanishing terms with exponents from the
isolated subgroup of rank 1. This does not take place for
valuations of type II.2 (exceptional curve valuations).

Each series from the ring $Z[[S]]$ with the free term
equal to $1$ has a unique representation of the form
$\prod\limits_{ a\in S}(1-t^{ a})^{k_{a}}$
with $k_a\in\Z$. (Generally speaking, this product is not
finite: infinite number of the exponents $k_a$ may be different
from zero.)

The Poincar\'e series of rank $1$ valuations have the following form.

1) For a divisorial valuation
\begin{equation}\label{Ps_div}
P(t)=\frac{\prod\limits_{i=1}^h(1-t^{m^{\alpha_i}})}
{\prod\limits_{i=0}^{h+1}(1-t^{m^{\beta_i}})}
\end{equation}
where the exponents $m^{\beta_i}$ generate an infinite cyclic
group (i.e. a group isomorphic to $\Z$). Here $\alpha_i$ are
rupture points and $\beta_i$ are dead ends of the dual graph of
the minimal resolution.

2) For a valuation of type II.4 (an irrational one) the
Poincar\'e series has the same form \ref{Ps_div}, but the
exponents $m^{\beta_i}$ generate a free abelian group of rank
2. (The ratios $m^{\beta_i}/m^{\beta_0}$ are rational for $i\le
h$ and the ratio $m^{\beta_{h+1}}/m^{\beta_0}$ is irrational.

3) For a valuation of type II.3 (an infinite valuation) the
Poincar\'e series has a representation of the form
$$
P(t)=
\frac{\prod\limits_{i=1}^{\infty}(1-t^{m^{\alpha_i}})}
{\prod\limits_{i=0}^{\infty}(1-t^{m^{\beta_i}})}
$$
with infinitely many factors.

This shows that rank $1$ valuations of different types cannot
have equal Poincar\'e series.

\begin{remark}
This also can be deduced from the fact that, for one valuation,
the set of exponents with non-zero coefficients coincides with
the semigroup of values of the valuation.
\end{remark}

For a collection of divisorial valuations the statement of the
Theorem was proved in \cite{FAOM}.
We shall reduce consideration of a collection of valuations
of  different types to the case of a collection of
divisorial valuations.
For that we shall substitute each non-divisorial valuation from
the collection by one or two divisorial ones. In each
case the  Poincar\'e series of the resulting collection of valuations
should be defined by the Poincar\'e series of the initial
collection and the dual graph of the minimal resolution of the
resulting collection (or rather os series of them) should permit
to restore
the dual graph of the minimal resolution of the initial one.

Suppose that a valuation from the collection $\{v_i\}$ is of type
II.1 (a curve valuation of rank $2$). Without loss of generality we
may assume that this is the first one. Let $E_{\sigma_0}$ be a
vertex of the dual graph of the minimal resolution of the collection
$\{v_i\}$ far enough on the infinite tail corresponding to the
valuation $v_{1}$. Let us substitute the (rank $2$) valuation
$v_{1}$ in the collection $\{v_i\}$ by the divisorial valuation
$v_{\sigma_0}$. If one knows the dual graph of the minimal
resolution of the new collection, one can easily restore the dual
graph for the initial one. One can see that the Poincar\'e series of
the new collection of valuations is obtained from the Poincar\'e
series $P_{\{v_i\}}(t_1, \ldots, t_s)$ (see equation~(\ref{eq2})) by
substituting $t_{1}^{(0,1)}$ by $t_{\sigma_0}$ and $t_{1}^{(1,0)}$
by $t_{\sigma_0}^N$ with $N$ large enough.

Suppose that the valuation $v_{1}$ is of type II.2 (an
exceptional curve valuation), let $E_{\sigma_0}$ be the
exceptional curve (a component of the exceptional divisor)
corresponding to the valuation, and let $P$ be the
corresponding point of the component $E_{\sigma_0}$. (The point
$P$ is the intersection point of the component $E_{\sigma_0}$
with another component of the exceptional divisor.) Let us make
sufficiently
many additional blowing-ups at the point $P$ of the component
$E_{\sigma_0}$. (If one knows the dual graph of a resolution
obtained this way, one can easily restore the
dual graph for the minimal one.) In the new resolution, let
$P$ be the intersection point of the component $E_{\sigma_0}$
with a component $E_{\sigma'_0}$. Let us substitute, in the
collection $\{v_i\}$, the valuation $v_{1}$ by two divisorial
valuations $v_{\sigma_0}$ and $v_{\sigma'_0}$. One can see
that the Poincar\'e series of the new collection of valuations
is obtained from the Poincar\'e series $P_{\{v_i\}}(t_1,
\ldots, t_s)$ by substituting $t_{1}^{(1,0)}$ by
$t_{\sigma_0}$ and $t_{1}^{(0,1)}$ by
$t_{\sigma'_0}t_{\sigma_0}^N$ with $N$ large enough.

Suppose that $v_1$ is a valuation of type I.3 (an
infinite one). Moreover let us assume that
$v_1,\seq{v}2r$ are all valuations of type I.3 in the collection.
The Poincar\'e series has the form
$$
\prod_{n=1}^{n_0} (1-\tt^{\,\kk^n})^{-\chi_n}
\prod_{\ell=1}^{r}
\prod_{j=1}^{\infty}
\frac{1-\tt^{\,\mm^{\alpha^{\ell}_{j}}}}
{1-\tt^{\,\mm^{\beta^{\ell}_{j}}}}
\; .
$$
Here the product
$\prod\limits_{j=1}^{\infty}
\frac{1-\tt^{\,\mm^{\alpha^{\ell}_{j}}}}
{1-\tt^{\,\mm^{\beta^{\ell}_{j}}}}$ corresponds to the infinite
part (tail) of the dual graph corresponding to the valuation
$v_{\ell}$. In these products one has
$$
\mm^{\alpha^{\ell}_{1}} <
\mm^{\beta^{\ell}_{1}} <
\mm^{\alpha^{\ell}_{2}} <
\mm^{\beta^{\ell}_{2}}   < \ldots
$$
and $\mm^{\alpha^{\ell}_{j}}\to \infty$, $\mm^{\beta^{\ell}_{j}}\to
\infty$ when $j\to \infty$. The distribution of the factors between
these products (from a certain place) is determined by the following
property: the ratio $m_{j_1}^{\beta^{\ell}_{j}}/
m_{j_2}^{\beta^{\ell}_{j}}$ strictly increases with $j$ along the
tail corresponding to the valuation $v_{j_1}$ (i.e. for $\ell=j_1$),
strictly decreases along the tail corresponding to the valuation
$v_{j_2}$ ($\ell=j_2$) and is constant along all other tails.

Let us substitute the infinite valuation $v_1$ by the
divisorial valuation $v_{\alpha^{1}_{N}}$ with $N$ large enough.
The dual graph of
the minimal resolution of the collections
$\{v_{\alpha^{1}_{N}}, \seq{v}2s\}$ is obtained
from the
dual graph of the minimal resolution of the collections $\{v_i\}$
by
truncation of the infinite tail corresponding to the valuation
$v_1$
at the vertex
$\alpha^{1}_{N}$.
The dual graph of the resolution of the collection
$\{v_i\}$ can be restored if one knows the truncated ones for all
$N$ large enough.
The Poincar\'e series of the collection
$\{v_{\alpha^{1}_{N}}, \seq{v}2s\}$
is equal to
$$
\prod_{n=1}^{n_0} (1-\tt^{\,\kk^n})^{-\chi_n}
\prod_{j=1}^{N}
\frac{1-\tt^{\,\mm^{\alpha^{1}_{j}}}}
{1-\tt^{\,\mm^{\beta^{1}_{j}}}}
\prod_{\ell=2}^{r}
\prod_{j=1}^{\infty}
\frac{1-\tt^{\,\mm^{\alpha^{\ell}_{j}}}}
{1-\tt^{\,\mm^{\beta^{\ell}_{j}}}}
\; .
$$

Suppose that the valuation $v_1$ is of type II.4 (an irrational
one). The semigroup of values of the valuation $v_1$ is generated by
$1 = m_1^{\beta_0}, m_1^{\beta_1},m_1^{\beta_2}, \ldots,
m_1^{\beta_h}$ where $m_1^{\beta_0}< m_1^{\beta_1} < m_1^{\beta_2} <
\ldots < m_1^{\beta_h}$, the numbers $m_1^{\beta_i}$ are rational
for $i<h$ and $m_1^{\beta_h}$ is irrational. Let us substitute the
valuation $v_1$ in the collections $\{v_i\}$ by the divisorial
valuation $v_{\sigma_h}$ far enough in the sequence $\{\sigma_h\}$
corresponding to the valuation $v_1$ and moreover such that in the
dual graph of the minimal resolution of the valuation $v_1$ the
component $E_{\sigma_h}$ does not intersect the component
$E_{\sigma_{h+1}}$ (the first among those born by blowing-ups at
points of the component $E_{\sigma_h}$). The dual graph of the
resolution of the collection $\{v_i\}$ can be restored if one knows
the dual graphs of the minimal resolutions of the collections
$\{v_{\sigma_h}, \seq{v}2s\}$ with all $h$ large enough of the
described type. The Poincar\'e series of the collection
$\{v_{\sigma_h},\seq v2s\}$ is obtained from the Poincar\'e series
$P_{\{v_i\}}(\tt)$ of the collection $\{v_i\}$ by substituting each
monomial $t_1^{m_1}t_2^{m_2}\ldots t_s^{m_s}$ by
$t_1^{m'_1}t_2^{m_2}\ldots t_s^{m_s}$ where $m'_1$ is defined in the
following way. Let $\alpha_N$ be the result of the truncation of the
continuous fraction of the (irrational) number $m_1^{\beta_h}$ at
the level $N$ large enough. Let $m_1 = k_0 m_1^{\beta_0} + \cdots +
k_{h-1}m_1^{\beta_{h-1}} + k_h m_1^{\beta_h}$. Then $m'_1 = k_0
m_1^{\beta_0} + \cdots + k_{h-1}m_1^{\beta_{h-1}} + k_h \alpha_N$.
\end{proof}

\end{document}